\documentstyle[12pt]{article}

\newfam\frakfam
\font\teneufm=eufm10     \textfont\frakfam=\teneufm
\font\seveneufm=eufm7    \scriptfont\frakfam=\seveneufm
\font\fiveeufm=eufm5     \scriptscriptfont\frakfam=\fiveeufm
\def\frak{\fam\frakfam \teneufm}

\newfam\msbfam
\font\tenmsb=msbm10      \textfont\msbfam=\tenmsb
\font\sevenmsb=msbm7     \scriptfont\msbfam=\sevenmsb
\font\fivemsb=msbm5      \scriptscriptfont\msbfam=\fivemsb
\def\Bbb{\fam\msbfam \tenmsb}

\font\titlebf=cmbx10  scaled \magstep2

\newcommand{\ie}{\mbox{\em i.e.}}
\newcommand{\fq}{\mbox{${\Bbb F}_q\,$}}
\newcommand{\fp}{\mbox{${\Bbb F}_p\,$}}
\newcommand{\clo}{\mbox{${\overline {\Bbb F}}_p\,$}}
\newcommand{\nq}{\mbox{${\Bbb F}_{q^n}\,$}}

\newcommand{\sn}{\mbox{${\frak S}_n\,$}}
\newcommand{\op}{\mbox{${\mathcal O}_P\,$}}

\newcommand{\pf}{\mbox{${\Bbb P}(F)\,$}}
\newcommand{\pe}{\mbox{${\Bbb P}(E)\,$}}

\newcommand{\jth}{\mbox{$j^{th}\;$}}
\newcommand{\ith}{\mbox{$i^{th}\;$}}

\newcommand{\degdiff}{\mbox{\rm degDiff}}
\newcommand{\gal}{\mbox{\rm Gal}}
\newcommand{\diff}{\mbox{$\frak D$}}

\newtheorem{theorem}{Theorem}[section]
\newtheorem{lemma}[theorem]{Lemma}
\newtheorem{proposition}[theorem]{Proposition}
\newtheorem{corollary}[theorem]{Corollary}
\newtheorem{definition}[theorem]{Definition}

\newtheorem{example}[theorem]{\sc Example}

\renewcommand{\labelenumi}{\rm (\roman{enumi})}

\begin{document}
\font\titlebf=cmbx10 scaled \magstep2
\centerline{\titlebf Symmetry, splitting rational places in extensions}
\centerline{\titlebf of function fields and generalization of the}
\centerline{\titlebf Hermitian function field}
\medskip
\centerline{{\bf Vinay Deolalikar} \footnote[1]{ This work was part of the author's Ph.D. thesis and was supervised by the late Prof. Dennis Estes.}}
\medskip
\medskip
\centerline{\sc Dedicated to the late Prof. Dennis Estes}
\centerline{Mathematics Subject Classification 14G05, 14G50}

\section{Introduction}

Let $F/K$ be an algebraic function field in one variable over a finite field of constants $K$, \ie, $F$ is a finite algebraic extension of $K(x)$ where $x \in F$ is transcendental over $K$. Let $E$ be a finite separable extension of $F$. Let $N(E)$ and $g(E)$ denote the number of places of degree one (or {\em rational} places), and the genus, respectively, of $E$. Let $[E:F]$ denote the degree of this extension.

In recent years, there has been a spurt of interest in  algebraic function fields with many rational places, or, equivalently, curves over finite fields with many rational points. The initial impetus for this interest came from applications to coding theory, wherein, in 1981, Goppa \cite{Gop1} discovered that function fields with  many rational places could be used to construct long codes, whose parameters could then be ascertained using the Riemann-Roch theorem. Since then, more applications of such function fields have been discovered \cite{Xin1}.  Equally importantly, function fields with many rational  places  are an interesting mathematical problem in their own right, with connections to several well studied problems in arithmetical algebraic geometry, and have been recognised as such. Consequently, various aspects of such function fields have been studied. Many authors have also written on this subject in the language of curves over finite fields.

One technique to produce function fields with many rational places is to somehow split many rational places of the projective line in an extension, while keeping the rise in genus low. The technique existing in literature that can be used to achieve this uses class field theory and was introduced by Serre \cite{Ser1, Ser2, Ser3}. For practical applications of such function fields, however, it is necessary that the constructions be explicit, in that generators and equations satisfied by them should be provided. In this paper, we construct infinite families of extensions of the projective line in which all finite rational places split completely, the rise in genus incurred is low, and moreover, the construction is explicit. Furthermore, the constructions described herein work not only over finite fields of square cardinality, but over any non-prime finite field. Our techniques also shed light on why the square cardinality constraint arises in many existing examples.

To these ends, we introduce the notion of ``$(n,q)$-symmetry'' which is related to the notion of symmetry in polynomial rings in several indeterminates. For $q$ a prime power and $n \geq 2$, there are $n$ $(n,q)$-symmetric polynomials, which include trace and norm as the first and $n^{th}$, respectively. We call these $(n,q)$-symmetric polynomials $\{s_{n,i}\}_{1 \leq i \leq n}$. These polynomials, among other things, have the property that they map \nq into \fq.  Using these symmetric polynomials, we then construct two sets (of $n-1$ families each) of function fields over ${\Bbb F}_{q^n}$, which we name $\{E_i\}_{2 \leq i \leq n}$ and $\{E_{j,Kum}\}_{1 \leq j \leq n-1}$. Function fields belonging to these families are called symmetric function fields. The well known trace-norm function field is merely the special case of the family $E_n$ of our construction. Indeed, we show that it is the {\em worst} of these families in the sense of $N/g$ ratios. Similarly, many of the existing examples of function fields with many rational places are explained in the language of $(n,q)$-symmetry, thus unifying several existing examples under one theory.

We show that the $N/g$ ratios for the $n-1$ symmetric families $\{E_i\}$ follow a monotone pattern, namely
$$ N(E_2) = N(E_3) = \ldots = N(E_n),$$
$$ g(E_2) > g(E_3) > \ldots > g(E_n).$$

We also show that using the techniques described in this paper, we can explicitly construct subcovers of the symmetric function fields, and these also lead to examples with good $N/g$ ratios. Accordingly, we construct two new function fields meeting the Oesterle bounds where one is a subcover of the other.

We also compare our constructions to existing class field theoretic constructions. We show that for identical patterns of ramification groups, symmetric function fields will have better $N/g$ ratios than the corresponding class field theoretic constructions. The function fields that we construct are not maximal - they sit inside of some Ray class fields. Furthermore, while most of the examples using Ray class fields (example, \cite{Lau2}) have degree greater than the cardinality of the underlying finite field, ours have degrees lesser than that cardinality. This ties in with the approach used by other authors, for example \cite{NieXin1}, where the minimal function field with a certain property is shown to have the best $N/g$ ratio.

Another point of considerable importance is that in many existing examples of explicitly constructed  extensions of function fields where most of the rational places split completely, the genus, and number of rational places, of the resulting function field are usually fairly low.  However, for applications to coding theory, one would like to have a large number of rational places in order to build long codes. The function field towers of Garcia and Stichtenoth \cite{GarSti1, GarSti2}, give us a structure to build long codes, but the problem of finding a basis for the vector spaces of regular functions on them is yet to be solved. The constructions provided in this paper can provide examples of extensions of function fields $E/F$ where all rational places except one split completely, for arbitrarily high degree of extension $[E:F]$, and therefore, for arbitrarily high values of $N(E)$.

Finally, we provide a ``generalization'' of the Hermitian function field over $K$ where $K$ is an odd degree extension  of ${\Bbb F}_q$. For a long time, it has been known that the Hermitian function field has many unique properties, such as its maximality in the Hasse-Weil sense, and its large automorphism group. However there has hitherto not been a satisfactory generalization of this function field to field of nonsquare cardinality. We propose a solution to this problem. We suggest that the trace-norm function field, which has sometimes been implicitly assumed to be the generalization of the Hermitian function field, is not the appropriate generalization. Rather it is the function field $E_2$ that generalizes the Hermitian.

The paper is organized as follows. In section~\ref{section:preliminaries}, we introduce the notation used in the rest of the paper. In section~\ref{section:symmetric-functions}, we introduce the notion of $(n,q)$-symmetry and explain how it is related to the notion of symmetry in polynomial rings in several indeterminates. We also derive some elementary properties of $(n,q)$-symmetric polynomials. In section~\ref{section:symmetric-extensions-of-function-fields}, we construct Artin-Schreier and Kummer extensions of the projective line over \nq using $(n,q)$-symmetric polynomials. We also prove the irreducibility of the defining  equations and derive the sequence of ramification groups in the corresponding extensions. These do not follow from basic Artin-Schreier theory since there is no rational place at which the valuation of the $1^{st}$ to the $(n-1)^{st}$ $(n,q)$-symmetric polynomials is coprime to the characterisitic. In section~\ref{section:trace-norm}, we show how the well-known trace-norm construction is merely a special case of our constructions. In section~\ref{section:rational-places-versus-genus}, we compare the $N/g$ ratios for the various symmetric families of function fields, and show that the trace-norm family has the lowest $N/g$ ratio, \ie, all the other families have higher $N/g$ ratios. In section~\ref{section:generalization-of-Hermitian}, we propose the function field $E_2$ as the generalization of the Hermitian function field over non-square finite fields. We do this by examining the genus, subgroups of the automorphism groups of the two function fields, and the number of places of degree greater than one.

\section{Preliminaries} \label{section:preliminaries}

Throughout this paper, we will use the following notation:

For symmetric polynomials:
\begin{tabbing}
\sn \hspace{1.6cm} \= the symmetric group on $n$ characters \\
${\bf s}_{n,i}(X)$ \> the \ith elementary symmetric polynomial on $n$ variables \\
$q$  \>  a power of a prime $p$ \\
${\Bbb F}_l$  \> the finite field of cardinality $l$ \\
$s_{n,i}(t)$ \> the \ith $(n,q)$-elementary symmetric polynomial
\end{tabbing}

For function fields and their symmetric subfields:
\begin{tabbing}
$K$ \hspace{1.6cm} \= the finite field of cardinality $q^n$, where $n>1$ \\
$F/K$ \> an algebraic function field in one variable whose full field of constants is $K$ \\
$F_s$ \> the subfield of $F$ comprising $(n,q)$-symmetric functions \\
$F_s^\phi$ \> the subfield of $F_s$ comprising functions whose coefficients are from \fq \\
$E$ \> a finite separable extension of $F$, $E=F(y)$ where $\varphi(y) = 0$
for some irreducible  \\
 \> polynomial $\varphi[T] \in F[T]$
\end{tabbing}

For a generic function field $F$:
\begin{tabbing}
${\Bbb P}(F)$  \hspace{1.3cm} \=  the set of places of $F$ \\
$N(F)$  \> the number of places of degree one in $F$ \\
$N_m(F)$\> the number of places of degree $m, m>1$, in $F$ \\
$g(K)$ \> the genus of $F$   \\
$P$ \> a generic place in $F$  \\
$v_P$ \> the normalized discrete valuation associated with the place $P$ \\
\op \> the valuation ring of the place $P$ \\
$P'$ \> a generic place lying above $P$ in a finite separable extension of $F$\\
$e(P'|P)$ \> the ramification index for $P'$ over $P$ \\
$d(P'|P)$ \> the different exponent for $P'$ over $P$
\end{tabbing}

For the rational function field $K(x)$:
\begin{tabbing}
$P_\alpha$ \hspace{1.6cm} \= the place in $K(x)$ that is the unique zero of $x-\alpha, \; \alpha \in K$ \\
$P_\infty$ \> the place in $K(x)$ that is the unique pole of $x$
\end{tabbing}

\section{Symmetric functions} \label{section:symmetric-functions}

Let $R$ be an integral domain and $\overline{R}$ its field of fractions. Consider the polynomial ring  in $n$ variables over $R$, given by  $R\,[X]= R\,[x_1,x_2,\ldots,x_n]$. The symmetric group \sn acts in a natural way on this ring by permuting the variables.

\begin{definition}
A polynomial ${\bf f}(X) \in R\,[X]$ is said to be symmetric if it is fixed under the action of \sn. If \sn is allowed to act on $\overline{R}(X)$ in the natural way, its fixed points will be called symmetric rational functions, or simply, symmetric functions. These form a subfield $\overline{R}(X)_s$  of $\overline{R}(x)$.
\end{definition}

Symmetric functions form a very elegant branch in the study of polynomials with several indeterminates. We recall here one of the classical results on symmetric functions, often known as {\em The fundamental theorem on symmetric functions}\footnote{A sharper result, called {\em The fundamental theorem on symmetric polynomials}, says that every symmetric polynomial can be written as a polynomial in the elementary symmetric functions \cite{BeaPie1}.} \cite{BeaPie1}.
\begin{theorem}[The fundamental theorem on symmetric functions]
$$\overline{R}(X)_s  = \overline{R}({\bf s}_{n,1}(X), {\bf s}_{n,2}(X), \ldots, {\bf s}_{n,n}(X)) $$
where $\{{\bf s}_{n,i}(X)\}_{1\leq i \leq n}$ are given below:
\begin{eqnarray*}
  {\bf s}_{n,1}(X) &=& \sum_{i=1}^{n} x_i, \\
  {\bf s}_{n,2}(X) &=& \sum_{i<j \atop 1 \leq i,j \leq n} x_ix_j, \\
     \vdots  & &  \vdots \\
  {\bf s}_{n,n}(X) &=& x_1x_2\ldots x_n.
\end{eqnarray*}
 ${\bf s}_{n,i}(X)$ is called the \ith elementary symmetric polynomial in $n$ variables.
\end{theorem}

We would like to apply these notions to the setting of polynomial rings in one variable over finite fields. We start with the following elementary proposition.
\begin{proposition} \label{proposition:extension-finite-fields}
Consider the extension of finite fields given by ${\Bbb F}_{q^n}/{\Bbb
F}_q$. Then, the following hold:
\begin{enumerate}
\item This is a Galois extension of degree $n$. The Galois group of this extension, $G=\gal({\Bbb F}_{q^n}/{\Bbb F}_q)$ is cyclic, and is generated by $\phi \in G, \; \phi:\alpha \rightarrow \alpha^q$. This generating element is also called the ``frobenius'' of this extension.
\item The fixed field of the subgroup of $G$ generated by $\phi^k$ is ${\Bbb F}_{q^{\gcd(k,n)}}$.
\end{enumerate}
\end{proposition}

\begin{definition}
For the extension ${\Bbb F}_{q^n}/{\Bbb F}_q$, we will evaluate the elementary symmetric polynomials (resp. symmetric functions) in ${\Bbb F}_{q^n}(X)$ at $(X)=(t,\phi(t),\ldots,\phi^{n-1}(t))=(t,t^q,\ldots,t^{q^{n-1}})$. These will be called the $(n,q)$-elementary symmetric polynomials (resp. $(n,q)$-symmetric functions). For ${\bf f}(X) \in \nq\!(X)$, we will denote ${\bf f}(t,t^q,\ldots,t^{q^{n-1}})$ by $f(t)$, or, when the context is clear, by $f$.
\end{definition}

Thus the  $(n,q)$-elementary symmetric polynomials are the following:
\begin{eqnarray*}
s_{n,1}(t) &=& \sum_{0 \leq i \leq n-1} t^{q^i}, \\
s_{n,2}(t) &=& \sum_{i<j \atop 0 \leq i,j \leq n-1}t^{q^i}t^{q^j}, \\
    \vdots & & \vdots  \\
s_{n,n}(t) &=& t^{1+q+q^2+\ldots+q^{n-1}}.
\end{eqnarray*}

As an example, here are the $5$ $(5,q)$-symmetric polynomials.
$$
\begin{array}{|c|c|} \hline
s_{5,1}\mbox{(trace)}& t^{q^4} + t^{q^3} + t^{q^2} + t^q + t \\ \hline
s_{5,2}        & t^{q^4+q^3} + t^{q^4 + q^2} + t^{q^4 + q} + t^{q^4 + 1} +
                 t^{q^3 + q^2} + t^{q^3 + q} \\
               & + t^{q^3 +1} + t^{q^2+q} + t^{q^2 +1} + t^{q+1} \\ \hline
s_{5,3}        & t^{q^4+q^3+q^2} + t^{q^4+q^3 +q}+t^{q^4 + q^3 + 1} + t^{q^4 +                  q^2+q} + t^{q^4 + q^2 + 1} \\
               & + t^{q^4 + q + 1} + t^{q^3+q^2+q} + t^{q^3+q^2 +1} + t^{q^3+q+1} + t^{q^2+q+1} \\ \hline
s_{5,4}        & t^{q^4+q^3+q^2+q} + t^{q^4+q^3+q^2+1} + t^{q^4+q^3+q+1} \\
               & + t^{q^4 + q^2 + q + 1} + t^{q^3+q^2+q+1} \\ \hline
s_{5,5}\mbox{(norm)}  & t^{q^4 + q^3+q^2+q+1} \\  \hline
\end{array}
$$

We now proceed to give some elementary properties of $(n,q)$-symmetric polynomials.
\begin{lemma} \label{lemma:fsphi-in-fq}
Let $f(t)$ be an $(n,q)$-symmetric function with coefficients from \fq and $\gamma \in \nq$. Then, we have that $f(\gamma) \in {\Bbb F}_q \cup \infty$.
\end{lemma}
{\bf Proof}. $f(t)$ can be written as a rational function in $\{s_{n,i}(t)\}_{1 \leq i \leq n}$, with coefficients from \fq.  By construction, each $(n,q)$-elementary symmetric polynomial is invariant under the operation of raising to the $q^{th}$ power, modulo $(t^{q^n} - t)$. In other words, each $(n,q)$-elementary symmetric polynomial, when restricted to \nq, is invariant under this operation. So are the coefficients, since they are chosen from \fq. The inclusion of infinity in the range comes since $\gamma$ may be a pole of $f$. \hfill $\Box$

There are $n$ $(n,q)$-elementary symmetric polynomials. $s_{n,1}$ and $s_{n,n}$  are also known as ``trace'' and ``norm,'' respectively. In the case of $n=2$, these are the only $(n,q)$-elementary symmetric polynomials. For $n \geq 3$, we have a greater choice, which has not been exploited in existing techniques to construct function fields with many rational places. The main thesis of this paper is that these additional polynomials are indeed the more useful ones, as $n$ increases. In particular, we show that a ``generalization'' of the Hermitian function field can be obtained for $n\geq3$ by using $s_{n,2}$.

We end this section with some lemmas that will be used later. All of the following are valid for $1 \leq i \leq n$.

\begin{lemma} \label{lemma:Sni-Sn(n-i)}
$$ \frac{s_{n,i}(t)}{s_{n,n}(t)} = s_{n,n-i}\left(\frac{1}{t}\right).$$
Thus, there is a bijection between the set of nonzero roots of $s_{n,i}(t)$ and that of $s_{n,n-i}(t)$. More precisely, $s_{n,i}(\alpha) = 0 \Leftrightarrow s_{n,n-i}(1/\alpha) = 0$, where $\alpha \neq 0$.
\end{lemma}
{\bf Proof}. Follows immediately from the relation between $s_{n,i}(t)$ and $s_{n,n-i}(t)$ and observing that the only root of $s_{n,n}(t)$ is zero itself.

\begin{lemma} \label{lemma:derivative-of-Sni}
$$[s_{n,i}(t)]' = {[s_{n-1,i-1}(t)]}^q. $$
\end{lemma}
{\bf Proof}. The only terms that will contribute to the derivative are those whose exponent is coprime to $q$, \ie, in which $1$ is a summand. There are ${n-1 \choose i-1}$ such terms and they are all distinct.  Moreover, each such term, divided by $t$, is just the $q^{th}$ power of a term in $s_{n-1,i-1}(t)$. It is easy to see that this correspondence is bijective by noting that there are exactly ${n-1 \choose i-1}$ terms in $s_{n-1,i-1}(t)$ as well. \hfill $\Box$

\begin{lemma} \label{lemma:Sni-splits-partially}
The roots of $s_{n,i}(t)$ of multiplicity coprime to $p$ lie in \nq.
\end{lemma}
{\bf Proof}. The proof rests on the observation that
\begin{eqnarray*}
 {s_{n,i}(t)}^q - s_{n,i}(t) &=& [t^{q^n} - t] {[s_{n-1,i-1}(t)]}^q, \\
                 &=& [t^{q^n} - t] {s_{n,i}(t)}', \\
\mbox{or,                  } s_{n,i}(t)[{s_{n,i}(t)}^{q-1} - 1] &=& [t^{q^n} - t] {s_{n,i}(t)}'.
\end{eqnarray*}
Now consider a root $\alpha$ of multiplicity $m$ of $s_{n,i}(t)$ in the algebraic closure of ${\Bbb F}_q$, where $\gcd(m,p)=1$. This is not a root of the factor $[{s_{n,i}(t)}^{q-1} - 1]$ and it can only occur to a multiplicity $m-1$ in ${s_{n,i}(t)}'$. So it must appear to multiplicity $1$ in the factor $[t^{q^n} - t]$ on the RHS. Thus it must lie in ${\Bbb F}_{q^n}$.  \hfill $\Box$

\begin{lemma}
Every root of $s_{n,i}(t)$ lies in a field ${\Bbb F}_{q^k}$, where $n-i+1 \leq k \leq n$.
\end{lemma}
{\bf Proof}. From Lemma~\ref{lemma:Sni-splits-partially}, we know that if a root of $s_{n,i}(t)$ is not also a root of $s_{n-1,i-1}(t)$, it must lie in \nq. If, on the other hand, it is a root of $s_{n-1,i-1}(t)$,  it either lies in ${\Bbb F}_{q^{n-1}}$ or is a root of $s_{n-2,i-2}$. Now we descend this way till we reach $s_{n-i+1,1}(t)$, all of whose roots lie in ${\Bbb F}_{q^{n-i+1}}$. Notice that $s_{n-i+1,1}(t)$ can have no multiple roots since its derivative is equal to $1$. \hfill $\Box$

\begin{lemma}
Let $\alpha$ be a root of $s_{n,i}(t)$ of multiplicity $m>1.$ Further let
$\gcd(m,p)=\gcd(m-1,p)=1$. Then $\alpha \in {\Bbb F}_q$.
\end{lemma}
{\bf Proof}.  From Lemma~\ref{lemma:derivative-of-Sni}, we know that the roots of the derivative of $s_{n,i}(t)$ are the same as those of $s_{n-1,i-1}(t)$. From Lemma~\ref{lemma:Sni-splits-partially}, we know  that under the hypothesis of the present lemma, the roots in question of $s_{n-1,i-1}(t)$ will lie in ${\Bbb F}_{q^{n-1}}$, while those of $s_{n,i}(t)$ will lie in \nq. For a common root then, it must be a root of both of these polynomials, and hence must lie in the intersection of both these fields, which is ${\Bbb F}_q$. \hfill $\Box$

\begin{corollary} \label{corollary:multiple-root-of-Sni}
If  $p$ does not divide both ${n \choose i}$ and ${n-1 \choose i-1}$, then
there are no non-zero roots of $s_{n,i}(t)$ of multiplicity $m>1$, such that
$\gcd(m,p)=\gcd(m-1,p)=1$.
\end{corollary}
{\bf Proof}. When restricted to ${\Bbb F}_q$, $s_{n,i}(t)={n \choose i}t^i$, and  $s_{n-1,i-1}(t)={n-1 \choose i-1}t^{i-1}$. Thus, if the root is non-zero, the only way both these polynomials can be zero is if the characteristic divides both ${n \choose i}$ and ${n-1 \choose i-1}$. \hfill $\Box$

\begin{lemma} \label{lemma:Sn1-is-permutation}
$s_{n,1}(t)$ is a permutation polynomial over ${\Bbb F}_{q^m}$ if $\gcd(m,n)=1$ and $p$ does not divide $n$.
\end{lemma}
{\bf Proof}. Suppose the pre-image under $s_{n,1}(t)$ of an element comprised of two distinct elements, then, since $s_{n,1}(t)$ is additive, their difference would be its root. But all the roots of $s_{n,1}(t)$ lie in \nq. If $p$ does not divide $n$, then none of these lie in \fq. Thus, under these hypotheses, $s_{n,1}(t)$ would not have any roots in ${\Bbb F}_{q^m}$, where $\gcd(m,n)=1$, and we would get a contradiction.\hfill $\Box$

Now we wish to use these notions in the setting of an algebraic function field in one variable $F/K$, where $K = \nq$, and $F$ is an algebraic extension of $K(x)$, where $x$ is transcendental over $K$.

\begin{definition} The field of $(n,q)$-symmetric functions with coefficients in \nq will be denoted
$$  F_s = \nq(s_{n,1}(x), s_{n,2}(x), \ldots, s_{n,n}(x)) \subset K(x),$$
and the field of $(n,q)$-symmetric functions with coefficients in \fq will be denoted
$$  F_s^\phi = \fq(s_{n,1}(x), s_{n,2}(x), \ldots, s_{n,n}(x)) \subset K(x),$$
where the superscript $\phi$ indicates that the values that these functions take on \nq are fixed by $\phi$. Thus, they lie in \fq, as per Lemma~$\ref{lemma:fsphi-in-fq}$.
\end{definition}

\section{Symmetric extensions of function fields} \label{section:symmetric-extensions-of-function-fields}

Let $F$ and $K$ be as described earlier. Let $E$ be a finite separable extension of $F$, generated by $y$, where $\varphi(y) = 0$, for $\varphi(T)$ an irreducible polynomial in $F[T]$.

In this section we will introduce families of extensions of $F$ whose generators satisfy explicit equations involving only $(n,q)$-symmetric functions. Let $y$ satisfy
                        $$ g(y) = f(x),$$
where $f, g \in F_s^\phi$. If $K=\nq$, this implies that in the residue field of a rational place,  although the class of $x$ and $y$ will assume  values in ${\Bbb F}_{q^n}\cup \infty$, that of $f(x)$ and $g(y)$ will assume values only in ${\Bbb F}_{q}\cup \infty$. Among the Galois extensions that such equations can produce are the two special cases of extensions of Artin-Schreier and Kummer type. In this paper, we will mainly investigate the case where $f(x)$ is an $(n,q)$-elementary symmetric polynomial. But the techniques to analyse the case of arbitrary $f \in  F_s^\phi$ remain the same.

\subsection{Symmetric extensions of Artin-Schreier type}

In classical Artin-Schreier extensions, the Galois group is a $1$-dimensional vector space over a subfield of the field of constants. We will, however, consider a modified Artin-Schreier extension, where the Galois group is isomorphic to the subgroup of the additive group of  ${\Bbb F}_{q^n}$ comprising elements whose trace in  \fq is zero.

Before we proceed to describe such modified Artin-Schreier extensions, we reproduce here some definitions and the relations between the different, ramification groups and the genus in extensions of function fields. References are \cite{ParSha1}, \cite{Ser4}, \cite{Sti1} and  \cite{Wei1}.

\begin{definition}
 Let $F'/F$ be a finite separable extension of function fields. Then, the different of $F'/F$, denoted $\diff(F'/F)$, is a divisor in $F'$ given by
$$\diff(F'/F) = \prod_{P' \in {\Bbb P}(F')}P'^{d(P'|P)},$$ where $d(P'|P)$ is the different exponent of the place $P'$ lying over $P \in \pf$. The degree of $\diff(F'/F)$, denoted $\degdiff(F'/F)$, is  the degree of this divisor.
\end{definition}

\begin{proposition}[Hurwitz-genus formula]
Let $F'/F$ be a finite separable extension of function fields. Then,
$$2g(F') - 2 = [F':F](2g(F) - 2) + \degdiff(F'/F).$$
\end{proposition}

\begin{proposition}[Transitivity of different]
Let $F \subseteq F' \subseteq F''$, where $F''/F'$ and $F'/F$ are both finite separable extensions. Let $P \subseteq P' \subseteq P''$, where $P, P'$ and $P''$ are places of $F, F'$ and $F''$, respectively. Then,
$$ \diff(F''/F) = \diff(F''/F') \diff(F'/F). $$
This yields the following relation between the different exponents
$$ d(P''|P) = d(P''|P') + e(P''|P')d(P'|P). $$
\end{proposition}

\begin{definition}[Ramification groups]
Let $F'/F$ be a Galois extension of function fields. Let $P \subseteq P'$, where  $P$ and $P'$ are places in $F$ and $F'$, respectively. The \ith ramification group of $G = \gal(F'/F)$ relative to $P'$ is $G_i = \{s \in G\; | \; s(v) \equiv v \bmod {P'}^{i+1}, \forall v \in {\mathcal O}_{P'} \}$. Then, $G_{-1}$ is the decomposition group, $G_0$ is the inertia group and $G_{-1}/G_0$ is $\gal(K'/K)$, where $K$ and $K'$ are the residue fields of $P$ and $P'$, respectively. Further, $G_i$ is a normal subgroup of $G_{-1}$ as well as $G_{i-1}$, for $i \geq 0$. $G_1$ is a $p$-group and $G_0/G_1$ is a cylic group of order coprime to $p$.
\end{definition}

\begin{proposition}[Hilbert's different formula]
Let $F'/F$ be a Galois extension of function fields. Let $P \subseteq P'$, where  $P$ and $P'$ are places in $F$ and $F'$, respectively. Then we have that
$$d(P'|P) = \sum_{i=0}^\infty |G_i| - 1. $$
\end{proposition}

\begin{proposition} \label{lemma:ramification-groups}
Let $F''/F$ be a Galois extension of function fields. Let $H$ be a normal subgroup of $G = \gal(F''/F)$, and $F'$ be its fixed field. Let $P, P'$ and $P''$ be places in $F$, $F'$ and $F''$, respectively, with $P \subset P' \subset P''$. Then the \ith ramification group of $\gal(F''/F')$ relative to $P''$ is $G_i\cap H$.
\end{proposition}

Now we proceed to describe extensions of modified Artin-Schreier type.

\begin{proposition} \label{proposition:Art-Sch}
Let $F/K$ be an algebraic function field, where $K=\nq$ is algebraically closed in $F$. Let $w \in F$ and assume that there exists a place $P \in \pf$ such that
$$ v_P(w) = -m, \, m > 0 \mbox{  and  } \gcd(m,q)  = 1.  $$
Then the polynomial $l(T)-w = a_{n-1}T^{q^{n-1}} + a_{n-2}T^{q^{n-2}}+ \ldots + a_0T - w \in F[T] $ is absolutely irreducible. Further, let $l(T)$ split into linear factors over $K$. Let $E=F(y)$ with
$$  a_{n-1}y^{q^{n-1}} + a_{n-2}y^{q^{n-2}}+ \ldots + a_0y = w. $$
Then the following hold:
\begin{enumerate}
\item $E/F$ is a Galois extension, with degree $[E:F] = q^{n-1}$. $\gal(E/F) = \{\sigma_\beta:y \rightarrow y + \beta\}_{l(\beta)=0}$.
\item $K$ is algebraically closed in $E$.
\item The place $P$ is totally ramified in $E$. Let the unique place of $E$ that lies above $P$ be $P'$. Then the different exponent $d(P'|P)$ in the extension $E/F$ is given by
$$ d(P'|P) = (q^{n-1}-1)(m+1).$$
\item Let $R \in \pf$, and $v_R(w) \geq 0$. Then $R$ is unramified in $E$.
\item If $a_{n-1}=\ldots=a_0=1$, and if $Q\in \pf$ is a zero of $w-\gamma$, with $\gamma \in \fq$. Then $Q$ splits completely in $E$.
\end{enumerate}
\end{proposition}
{\bf Proof}. For (i) - (iv), pl. refer \cite{Sul1}. For (v), notice that under the hypotheses, the equation $ T^{q^{n-1}} +T^{q^{n-2}}+ \ldots + T = \gamma $ has $q^{n-1}$ distinct roots in $K$.  \hfill $\Box$

For many of the extensions that we will describe, there exists no place where the hypothesis of Proposition~\ref{proposition:Art-Sch} is satisfied, namely, that the valuation of $w$ at the place is negative and coprime to the characteristic. In particular, we need a criterion for determining the irreducibility of the equations that we will need to use. We provide such a criterion in Proposition~\ref{proposition:irreducibility} and Corollary~\ref{corollary:irreducibility}.

\begin{proposition} \label{proposition:irreducibility}
Let $V$ be a finite subgroup of the additive group of \clo. Then $V$ is a \fp-vector space. Define $L_V(T) = \prod_{v \in V}(T-v)$. Thus, $L_V(T)$ is a separable \fp-linear polynomial whose degree is the cardinality of $V$. Now let $h(T,x) =  L_V(T) - f(x)$, where $f(x) \in \clo[x]$. Then, $h(T,x)=L_V(T) - f(x)$ is reducible over $\clo[T,x]$ iff there exists a polynomial $g(x) \in \clo[x]$ and a proper additive subgroup $W$ of $V$ such that $f(x) = L_{W'}(g(x))$, where $W' = L_W(V)$.
\end{proposition}
{\bf Proof}. Suppose $h(T,x) = g_1(T,x)\ldots g_n(T,x)$ is a factorization of
$h(T,x)$ over $\clo[T,x]$. From Gauss' lemma, and the fact that $h(T,x)$ is monic in $T$, we see that this is equivalent to a factorization into monic irreducible factors in $\clo(x)[T]$.

Let $F = \clo(x)$ and view $h(T)$ and $g_i(T)$ as polynomials over $F$. Let $L = F[T]/(g_1(T))$. Then $L$ is a field. Let $y$ be the corresponding root of $g_1$ in $L$.

Note that $v \in V$ acts on the roots of $h(T)$ over $L$ via $(T \rightarrow T + v)$. Since $y$ is a root of $h(T)$, we see that $h(T)$ splits completely over $L$, into factors $(T-y-v)$, where $v \in V$.

Also, since $h(T)$ is separable, the factors $g_1(T),\ldots,g_n(T)$ are all distinct. Since the roots of $g_i(T)$ are all conjugate over $F$, it follows that the action of $V$ on the roots of $h(T)$ must respect the $g_i$. Thus, the action of $V$ induces an action on the set $\{g_1(T),\ldots,g_n(T)\}$.

Let $W$ be the stabilizer of $g_1$ in $V$. Since $V$ is transitive on $\{g_1(T),\ldots,g_n(T)\}$, $|W| = |V/n| = deg(g_1)$. Thus $W$ is isomorphic to the Galois group of $g_1(T)$. It follows that
\begin{eqnarray*}
g_1(T) & = & \prod_{w \in W} (T-y-w), \\
       & = & L_W(T-y), \\
       & = & L_W(T) - L_W(y).
\end{eqnarray*}
Also note that the factorization of $h(T)$ is nontrivial iff $W$ is a proper
subgroup of $V$.

Now note that the constant term in $g_1(T)$ is $L_W(y)$. Thus $L_W(y) \in F$. In particular, $L_W(y)$ is a polynomial in $x$ with \clo coefficients.  It is useful to note that $L_W : \clo \rightarrow \clo$ is a homomorphism of additive groups with kernel $W$. Finally, let $\hat{W}$ be a complement of $W$ in $V$ and let $W' = L_W(V)$. Then
\begin{eqnarray*}
h(T) & = & \prod_{v \in \hat{W}}\prod_{w \in W} (T - y - v - w), \\
     & = & \prod_{v \in \hat{W}} L_W(T - y - v), \\
     & = & \prod_{v \in \hat{W}} L_W(T) - L_W(y) - L_W(v), \\
     & = & \prod_{w \in W'}  L_W(T) - L_W(y) - w, \\
     & = & L_{W'}(L_W(T) - L_W(y)), \\
     & = & L_{W'}(L_W(T)) -L_{W'}(L_W(y)), \\
     & = & L_V(T) - L_{W'}(L_W(y)).
\end{eqnarray*}

Comparing this to the equation $h(T) = L_V(T) - f$, we see that $f = L_{W'}(L_W(y))$. This proves the proposition. \hfill $\Box$

The following observation also follows from the proof of Proposition~\ref{proposition:irreducibility}. Let $E = F(y)$ be an extension of $F = \clo(x)$ with Galois group $V$, such that $y$ satisfies the equation $L_V(y) = f(x)$. Let $W$ be a subgroup of index $p$ in $V$, and let $M$ be its fixed field. Then clearly, $z = L_W(y) \in M$. Let $W' = L_W(V) \cong V/W$. Then we have that
$$ L_{W'}(z) = L_{W'}(L_W(y)) = L_V(y) = f(x). $$
By noting that the previous equation is of degree $p$, it follows that $M = F(z)$, where $z$ satisfies the equation
$$ L_{W'}(z) = f(x). $$

\begin{definition}
For $f(x) \in \clo[x]$, a coprime term of $f$ is a term with non-zero coefficient in $f$ whose degree is coprime to $p$. The coprime degree of $f$ is the degree of the coprime term of $f$ having the largest degree.
\end{definition}

\begin{corollary} \label{corollary:irreducibility}
Let $f(x) \in \clo[x]$. Let there be a coprime term in $f(x)$ of degree $d$, such that there are no terms of degree $dp^i$ for $i>0$ in $f(x)$. Then $L_V(T) -f(x)$ is irreducible for any subgroup $V \subset \clo$.
\end{corollary}
{\bf Proof}. Suppose $f(x)$ is the image of a linear polynomial $\sum a_nx^{p_n}$. Then the coprime term can only occur in the image of the term $a_0x$. But then, the images of the coprime term under $a_nx^{p^n}$, for $n>0$ will have degrees that contradict the hypothesis.

\begin{example}
The equation
$$ y^{q^2} + y^q + y = x^{q^2 + q} + x^{q^2 + 1} + x^{q+1} $$
is absolutely irreducible over \clo. This follows from Corollary~$\ref{corollary:irreducibility}$ by noting that the coprime degree of the RHS is $q^2+1$, and there are no terms of degree $(q^2 + 1)p^i$, for $i>0$.
\end{example}

To state the main theorem of this section which describes a general extension of modified Artin-Schreier type using $(n,q)$-elementary symmetric functions, we need  the following lemma.

\begin{lemma} \label{lemma:subextensions}
 Let $F = K(x)$, where $K=\nq, q = p^m, r=m(n-1),$ and $E = F(y)$, where $y$ satisifes the following equation:
$$ y^{q^{n-1}} + y^{q^{n-2}} + \ldots + y = f(x), $$
and $f(x) \in F$ is not the image of any element in $F$ under a linear polynomial.
Then the following hold:
\begin{enumerate}
\item $E/F$ is a Galois extension of degree $[E:F]=q^{n-1}$. $\gal(E/F) = \{\sigma_\beta: y \rightarrow y + \beta \}_{s_{n,1}(\beta)=0}$ can be identified with the set of elements in \nq whose trace in \fq is zero by $\sigma_\beta \leftrightarrow \beta$. This gives it the structure of a $r$-dimensional \fp vector space.
\item There exists a tower of subextensions
$$ F=E^0 \subset E^1 \subset \ldots \subset E^{r} = E, $$
such that for $0 \leq i \leq r-1,\; [E^{i+1}: E^i]$ is a Galois extension of degree $p$.
\item Let $\{b_i\}_{1 \leq i \leq r}$ be a \fp-basis for $\gal(E/F)$. Then we can build the tower of subextensions as follows. We set $E^j$ to be  the fixed field of the subgroup of the Galois group that corresponds to the \fp-subspace generated by $\{b_1,b_2,\ldots, b_{r-j}\}$. Then, the generators of $E^{j}$ are $\{y_1,y_2,\ldots,y_j\}$, where $y_1,y_2,\ldots,y_{r}=y$ satisfy the following relations:
\begin{eqnarray*}
       y^p - B_{r}^{p-1} y &=& y_{r-1}, \\
       y_{r-1}^p - B_{r-1}^{p-1} y_{r-1} &=& y_{r-2}, \\
          \vdots &  & \vdots \\
       y_1^p - B_{1}^{p-1}y_1 &=& f(x),
\end{eqnarray*}
where,
$$
\begin{array}{rcll}
  \beta_{r,j} &=& b_{r-j+1},  & B_{r} = \beta_{r,r}, \\
  \beta_{r-1,j} &=& \beta_{r,j}^p - B_{r}^{p-1}\beta_{r,j}, & B_{r-1} = \beta_{r-1,r-1},\\
 \vdots & & \vdots &\vdots  \\
  \beta_{1,j}   &=& \beta_{2,j}^p - B_2^{p-1}\beta_{2,j}, &  B_1 = \beta_{1,1}.
\end{array}
$$
\end{enumerate}
\end{lemma}
{\bf Proof}. For (i) refer Proposition~\ref{proposition:Art-Sch} and Proposition~\ref{proposition:irreducibility}. For (ii), note that since  $\gal(E/F)$ is an elementary Abelian group of exponent $p$, we can always find a normal series $\gal(E/K)=G^0 \rhd G^1 \rhd \ldots \rhd G^{r} = 1$ such that $|G^{i+1}/G^i| = p$. Now the existence of a tower $E^0 \subset E^1 \subset \ldots \subset E^{r}$, with $E^{j+1}/E^j$ Galois of degree $p$ is guaranteed by the Fundamental Theorem of Galois Theory by setting $E^i$ to be the fixed field of $G^{i}$. Alternately, we can get a constructive proof for (ii) from the proof for (iii), which follows.

For (iii), we set $G^{r-1}$ to be the \fp-span of $b_1$ and $E^{r-1}$ to be its fixed field. Then, the automorphisms of $E^{r}/E^{r-1}$ are given by $y \rightarrow y + ab_1$, with $a \in \fp$. Thus we have that
$$ \prod_{a \in \fp} y - ab_1 = {b_1}^p \prod_{a \in \fp} \frac{y}{b_1} - a  = y^p - {b_1}^{p-1}y = y_{r-1}. $$
Now, we set $G^{r-2}$ to be the \fp-span of $\{b_1, b_2\}$, and iterate this procedure for $E^{r-1}/E^{r-2}$, keeping in mind that the automorphism of $E/F$ given by $y \rightarrow y + b_2$, when pulled down to an automorphism of $E^{r-1}/E^{r-2}$, is given by $y_{r-1} \rightarrow y_{r-1} + {b_2}^p - {b_1}^{p-1}b_2$. By similarly setting $G^{r-j}, j \geq 2$ to be the \fp-span of $\{b_1,b_2,\ldots,b_{j}\}$, and  pulling  down the appropriate automorphisms of $E/F$ to those of  $E^{r-j+1}/E^{r-j}$ we get the other defining equations. The observation made after the proof of Proposition~\ref{proposition:irreducibility} completes the proof. \hfill $\Box$

\begin{corollary} \label{corollary:p-extension}
Let $E,F,K$ be as in Lemma~\ref{lemma:subextensions}. Then, every subextension $E^1$ of $E$ which has degree $p$ over $F$ is of the form $F(z)$, where $z$ satisfies an equation
$$ z^p - Az = f(x),  $$
where $A \in \nq$.
\end{corollary}
It is useful to note that since $E^1$ is the fixed field of a subgroup of $\gal(E/F)$ of index $p$, we can obtain a basis for $\gal(E/F)$ by adding one more element to a basis for this subgroup.

We are now in a position to state our main theorem for this section.

\begin{theorem} \label{theorem:Art-Sch-symmetric}
Let $F=K(x)$ where $K=\nq, q =p^m \mbox{ and } r=m(n-1)$. Consider the family of extensions $E_i=F(y)$, $2 \leq i \leq n$, where $y$ satisfies the equation
\begin{equation}
y^{q^{n-1}} + y^{q^{n-2}} + \ldots+ y = s_{n,i}(x). \label{equation:Art-Sch-symmetric}
\end{equation}
Then the following hold:
\begin{enumerate}
\item $E_i/F$ is a Galois extension, with degree $[E_i:F] = q^{n-1}$. $\gal(E_i/F)= \{\sigma_\beta:y \rightarrow y + \beta\}_{s_{n,1}(\beta) = 0}.$
\item The only place of $F$ that is ramified  in $E_i$ is the unique pole  $P_{\infty}$ of $x$. Furthermore, $P_\infty$ is totally ramified in $E_i$. Let $P'_\infty$ denote the unique place of ${\Bbb P}(E_i)$ that lies above $P_\infty$.
\item  Let $m_i$ denote the coprime degree of $s_{n,i}(x)$. We have that
$$ m_i = q^{n-1} + q^{n-2} + \ldots + q^{n-i+1} + 1. $$
The filtration of ramification groups relative to $P_\infty'$ is as follows:
$$ \gal(E_i/F) = G_0 = G_1 = \ldots = G_{m_i+1}, $$
$$ G_{m_i + 2} = \{0\}.$$
\item The different exponent $d(P'_\infty|P_\infty)$ is given by
$$ d(P_\infty|P'_\infty) = (q^{n-1}-1)(m_i+1).$$
\item The genus of $E_i$ is given by
$$ g(E_i) = \frac{(q^{n-1}-1)(m_i-1)}{2}.$$
\item All other rational places of $F$ split completely in $E_i$ giving
 $$ N(E_i) = q^{2n-1}+1.$$
Thus,  the number of rational places is independent of the choice of $(n,q)$-elementary symmetric polynomial of $x$.
\end{enumerate}
\end{theorem}
{\bf Proof}. (i) follows from Corollary~\ref{corollary:irreducibility} by observing that $p$ times the coprime degree of $s_{n,i}(x)$ is greater than its degree. Let $G = \gal(E_i/F)$. To determine the sequence of ramification groups of $G$ relative to $P_\infty'$, we will use the function $i_G$ defined on the elements of the $G$ as follows:
$$ i_G(s) \geq k+1 \Leftrightarrow s \in G_k.$$

For every subgroup $H$ of $G$, define $(G/H)_k$ to be the $k^{th}$ ramification group of the fixed field of $H$, relative to its unique place lying above $P_\infty$. Then, define the function $i_{G/H}$ on the elements of $G/H$ as follows:
$$ i_{G/H}(\overline{s}) \geq k+1 \Leftrightarrow \overline{s} \in (G/H)_k. $$

Note that for all subgroups $H$ of index $p$ in $G$, we know that
\begin{eqnarray*}
i_{G/H}(\overline{s}) &=& \left\{ \begin{array}{ll} \infty & \mbox{if } \overline{s} = 0, \\
                    m_i + 2 & \mbox{else}.
\end{array}
\right.
\end{eqnarray*}

We claim that if $K$ is any proper subgroup of $G$ then the average of the values of the function $i_{G/K}$ on the non-zero elements of $G/K$,
$$ \displaystyle{\mbox{avg}}_{s\in G/K, s \neq 0} i_{G/K}(s) = m_i + 2. $$

To see this, suppose that $|G/K| = p^l$. Consider the subgroups $H$ such that $K \subset H \subset G$, with $[G:H]=p$. These are in $1:1$ correspondence with the subgroups of index $p$ in $G/K$, which number $\frac{p^l - 1}{p-1}$. A non-zero element $s \in G/K$ is contained in exactly $\frac{p^{l-1}-1}{p-1}$ of these. Now, from \cite{Ser4}, Ch. IV, Proposition 3, we get
$$ i_{G/H}(\overline{s}) = \frac{1}{p^{l-1}} \sum_{s \rightarrow \overline{s}} i_{G/K}(s). $$
Since each non-zero $s \in G/K$ is nontrivial in exactly $p^{l-1}$ of the $G/H$, so that summing over all non-zero $s \in G/K$, we get
$$ \sum_{K\subset H \subset G, [G:H]=p \atop \overline{s} \in G/H, \overline{s} \neq 0} i_{G/H}(\overline{s}) = \sum _{s\in G/K} i_{G/K}(s). $$

But as noted previously, each $i_{G/H}(\overline{s}) = m_i +2$. Since each side in the previous equation has the same number of terms, the average of the RHS = average of the LHS = $m_i + 2$. This proves the claim.

It now follows that $i_G(s) = m_i + 2$ for all non-zero $s \in G$. For
suppose not. Then there must exist $s \in G$ such that $i_G(s) > m_i + 2$. Since $i_G$ is constant on cyclic subgroups, we see that $\langle{s}\rangle \neq G$. But then the average on $G/\langle{s}\rangle$ will be less than $m_i + 2$, giving a contradiction.

(iv) now follows from Hilbert's different formula. (v) is a straightforward application of the Hurwitz-genus formula. (vi) follows from Proposition~\ref{proposition:Art-Sch}. We obtain the number of rational places in $E_i$ as follows. Each of the $q^n$ finite places in $F$ splits completely in $E_i$, giving $(q^n)(q^{n-1}) = q^{2n-1}$ rational places in $E_i$. Also, $P_\infty$ ramifies totally in $E_i$, giving a sum total of $q^{2n-1}+1$ rational places in $E_i$. \hfill $\Box$

\begin{example} \label{example:n=3,i=2}\rm (n=3, i=2) Let $F=K(x)$, where $K = {\Bbb F}_{q^3}$. Let $E=F(y)$ where $y$ satisfies the equation
$$ y^{q^2} + y^q + y = x^{q^2 + q} + x^{q^2 + 1} + x^{q+1}.$$
All rational places of $F$, except $P_\infty$,  split completely in $E$, giving a total of $(q^3)(q^2) = q^5$ rational places. Also, $P_\infty$ ramifies totally in $E$ giving one rational place. Thus, $N(E) = q^5+1$ rational places. The genus $g(E) =  \frac{(q^2-1)(q^2)}{2}$. In a later section we will see that this extension attains the Oesterle lower bound on genus for $q=2$, \ie, over ${\Bbb F}_{q^3} = {\Bbb F}_8$.
\end{example}

\begin{example} \label{example:n=3,i=2,subfield}\rm (n=3, i=2) In this example, we discuss a subfield of the Example~\ref{example:n=3,i=2}. Let $E,F,K$ all be as in Example~\ref{example:n=3,i=2}. Then consider the extension $E^1/F$ where $E^1=F(y_1)$, and $y_1$ satisfies the equation
$$ y_1^q + (1+b^{q^2-q})y_1 = x^{q^2 + q} +  x^{q^2 + 1} + x^{q+1},$$
where $b(\neq 0) \in K$ is an element whose trace in \fq is zero. Then $E^1$ is a subfield of $E$ and $E = E^1(y)$, with $y$ satisfying
$$ y^q - b^{q-1}y = y_1. $$
$P_\infty$ is totally ramified in $E^1$ and all other rational places of $F$ split completely in $E^1$ (This is clear since $E^1$ is a subfield of $E$, in which these statements are true). Let $P^1_{\infty}$ be the unique place of $E^1$ lying above $P_\infty$. Then we have that $d(P^1_{\infty}|P) = (q-1)(q^2 + 2) $ and $g(E^1) = \frac{(q-1)(q^2)}{2}$. Also, $N(E^1) = q^4+1$. $E^1$ also attains the Oesterle lower bound on genus for $q=2$.
\end{example}

\begin{example} \rm (n=4, i=3) Let $F=K(x)$, where $K = {\Bbb F}_{q^4}$. Let  $E=F(y)$ where $y$ satisfies the equation
$$ y^{q^3} + y^{q^2} + y^q + y = x^{q^3 + q^2 + q} + x^{q^3 + q^2 + 1} +
x^{q^3+q+1} + x^{q^2 + q + 1}. $$
All rational places of $F$, except $P_\infty$,  split completely in $E$. $P_\infty$ ramifies totally in $E$. This gives us $N(E)=q^7+1$. The genus $g(E) =  \frac{(q^3-1)(q^3+q^2)}{2}$.
\end{example}

\begin{definition}
A function field $F/K$, where $K = \nq$ is said to be median if $N(F) = q^n + 1$. Thus, the number of its rational places is exactly in the middle of the range allowed by the Hasse-Weil bound.
\end{definition}

\begin{proposition}
Let $F=K(x)$, where $K = {\Bbb F}_{q^m}, \; \gcd(m,n) = \gcd(p,n) =1$. Let $E_i = F(y)$ where $y$ satisfies the equation
$$ y^{q^{n-1}} + y^{q^{n-2}} + \ldots + y = s_{n,i}(x). $$
Then $E_i$ is median\footnote{Function fields which are median over infinitely many extensions of their field of constants (\fq in our case) are called {\em exceptional}. There is a rich theory to such function fields. For instance, it is known \cite{Fri1} that the roots of their zeta-function occur in cliques as roots of unity times each other.}. It retains this property if we replace $s_{n,i}(x)$ with any other polynomial such that the resulting equation is irreducible.
\end{proposition}
{\bf Proof}. From Lemma~\ref{lemma:Sn1-is-permutation}, we know that under the hypothesis,  $ y^{q^{n-1}} + y^{q^{n-2}} + \ldots+ y$ simply permutes the elements of ${\Bbb F}_{q^m}$.  Then, for every value that $x$ can take in ${\Bbb F}_{q^m}$, we have exactly one solution for $y$.  \hfill $\Box$

\begin{example}\rm The following are the $4$ symmetric Artin-Schreier extensions of ${\Bbb F}_{q^5}(x)$. Notice that
$E_5$ is the familiar trace-norm construction. All of these have the same number of rational points, and the
genus increases in the order
$$ g(E_2)< g(E_3) < g(E_4) < g(E_5).$$
In other words, the trace-norm has the worst (lowest) $N/g$ ratio.
$$
\begin{array}{|c|c|} \hline
E_2 & y^{q^4} + y^{q^3} + y^{q^2} + y^{q} + 1 = x^{q^4+q^3} + x^{q^4 + q^2} + x^{q^4 + q} + x^{q^4 + 1} + \\
    & x^{q^3 + q^2} + x^{q^3 + q} + x^{q^3 +1} + x^{q^2+q} + x^{q^2 +1} + x^{q+1} \\ \hline
E_3 & y^{q^4} + y^{q^3} + y^{q^2} + y^{q} + 1 = x^{q^4+q^3+q^2} + x^{q^4+q^3 +q}+x^{q^4 + q^3 + 1} + x^{q^4 + q^2+q} \\
    & + x^{q^4 + q^2 + 1} + x^{q^4 + q + 1} + x^{q^3+q^2+q} + x^{q^3+q^2 +1} + x^{q^3+q+1} + x^{q^2+q+1} \\ \hline
E_4 & y^{q^4} + y^{q^3} + y^{q^2} + y^{q} + 1 = x^{q^4+q^3+q^2+q} + x^{q^4+q^3+q^2+1} \\
    & + x^{q^4+q^3+q+1} + x^{q^4 + q^2 + q + 1} + x^{q^3+q^2+q+1} \\ \hline
E_5 & \mbox{(trace-norm)} y^{q^4} + y^{q^3} + y^{q^2} + y^{q} + 1 = x^{q^4 + q^3+q^2+q+1} \\  \hline
\end{array}
$$
\end{example}

\subsubsection{Special case: $i=n$ (trace-norm)} \label{section:trace-norm}

It has been known that in extensions of the form $E/F$ where $E=F(y)$ with
$$  s_{n,1}(y) = s_{n,n}(x),  $$
all the rational places of $F$, except $P_\infty$, split completely, yielding many rational places in $E$. Extensions of this form have been referred to as ``trace-norm'' extensions. However, we will treat this type of extension as a special case of the generalized symmetric extensions of Theorem~\ref{theorem:Art-Sch-symmetric}. In the notation of Theorem~\ref{theorem:Art-Sch-symmetric}, the trace-norm extension is $E_n$.

The most famous example of a trace-norm construction is the Hermitian function field with full field of constants ${\Bbb F}_{q^2}$ (\ie, the case $n = 2$), where the trace and norm are taken down to ${\Bbb F}_q$. This is a maximal function field in the Hasse-Weil sense. Whenever $q \neq p$, we may take the trace and norm down to a smaller subfield of ${\Bbb F}_{q}$. All three - the degree of the extension, the number of rational places, and the genus - for such a construction increase as we take traces and norms to smaller  subfields. The maximum for each of these is attained when we take the trace and norm down to the prime field ${\Bbb F}_p$.

For a function field over $K$, we may then construct trace-norm extensions by taking these into any subfield of $K$. In the language of $(n/m,q^m)$-elementary symmetric polynomials, the most general form of the trace-norm extension is given by
$$ s_{\frac{n}{m},1}(y) = s_{\frac{n}{m},\frac{n}{m}}(x), $$
where the trace and norm are being taken in the subfield ${\Bbb F}_{q^m}$.
Now we would like to see how the ratio $N/g$ in such extensions varies as we vary $m$ from $1$ to the value of the greatest proper divisor of $n$. In other words, we would like to see how this ratio varies as we change the subfield of $K$ in which we take the trace and norm.

\begin{lemma} \label{lemma:N/g-trace-norm}
Let $F=K(x)$ where $K=\nq$. Let $m \neq n, m|n$ and $\frac{n}{m} = r$. Let $E=F(y)$, where $y$ satisfies the equation
\begin{equation}
y^{q^{m(r-1)}} + y^{q^{m(r-2)}} + \ldots + y = x^{q^{m(r-1)}+ \ldots + q^m + 1}. \label{equation:N/g-trace-norm}
\end{equation}
Thus, the trace and norm are being taken to the field ${\Bbb F}_{q^m}$. Then the following holds:
\begin{enumerate}
\item The ratio $N/g$, of the number of rational places to genus, decreases with decreasing $m$.
\item If $n \equiv 0\bmod 2$, the maximum value of  $N/g$ for the extension of this type is obtained by taking trace and norm down to the field of cardinality $q^{\frac{n}{2}}$.
\end{enumerate}
\end{lemma}
{\bf Proof}. In the general case, $N(E) = q^{2n-m}+1$, and $g(E) = \frac{q^m(q^{n-m} - 1)^2}{2(q^m-1)}$. Thus,
\begin{eqnarray*}
 \frac{N(E)}{2g(E)} & = &  \frac{(q^m-1)(q^m+q^{2n})}{(q^n - q^m)^2}.
\end{eqnarray*}
The numerator increases with increasing $m$ while the denominator decreases with increasing $m$. The result follows.   \hfill $\Box$

\noindent{{\bf Note}: Lemma~\ref{lemma:N/g-trace-norm}, (ii) is the case of the Hermitian function field.}

\begin{corollary}
For extensions of the form given by Lemma$~\ref{lemma:N/g-trace-norm}$, the lowest value of the ratio $N/g$ is reached when we take norms and traces down to ${\Bbb F}_p$.
\end{corollary}

There are other extensions of function fields using trace and norm, and one of them is given below.

\begin{example} \rm Let $F=K(x)$, where $K={\Bbb F}_{q^2}$. Let $E=F(y)$ with
$$ y^q + y = \frac{x^{q+1}}{x^q + x}. $$
This is the function field at the second step of the tower of function fields attaining the Drinfeld-Vladut bound from \cite{GarSti2}. The only places of \pf that are ramified in $E$ are $P_{\infty}$ and $\{P_{\alpha}\}_{s_{2,1}(\alpha) = 0, \alpha \neq 0}$. These are totally ramified, each with different exponent $2(q-1)$. All other rational places split completely. Thus we have that $N(E) = q^3-q^2+2q$ and $g(E) = (q-1)^2$.
\end{example}

\begin{lemma}
The extension $E_n$ described in Theorem~$\ref{theorem:Art-Sch-symmetric}$ attains the Hasse-Weil bound for $n=2$, for all values of $q$. For $n>2$, it does not attain the Hasse-Weil bound for any value of $q$.
\end{lemma}
{\bf Proof}. Observe that for $n>2$,
$$  \frac{q{(q^{n-1} - 1)}^2}{2(q-1)} > \frac{q^{\frac{n}{2}}(q^{\frac{n}{2}} -
 1)}{2}. $$
The lemma then follows from a well-known result that says that a function field over ${\Bbb F}_l$ cannot attain the Hasse-Weil bound for genus $g>\frac{\sqrt{l}(\sqrt{l}-1)}{2}$, cf. \cite{Sti1}, Ch. V.3.3. \hfill $\Box$

\subsection{Symmetric extensions of Kummer type}
We now study extensions whose Galois group is a subgroup of the
multiplicative group $K^*$. For this we will need that the field
contain a primitive \jth root of unity $\xi_j$ for some $j$ coprime to $p$. In particular we know that $K$ contains $\xi_j$ for $j = \frac{q^n-1}{q-1}$.

\begin{theorem} \label{theorem:Kummer-symmetric}
Let $F=K(x)$ where $K=\nq$. Let  $E_{i,Kum}=F(y)$, $1\leq i \leq n-1$, where $y$ satisfies the equation
\begin{equation}
 y^{\frac{q^n-1}{q-1}} = s_{n,i}(x). \label{equation:Kummer-symmetric}
\end{equation}
Then the following hold:
\begin{enumerate}
\item $E_{i,Kum}/F$ is a cyclic Galois extension, with degree $[E_{i,Kum}:F] = \frac{q^n-1}{q-1}$. $\gal(E_{i,Kum}/F) = \{ \sigma_j: y \rightarrow y{\xi}^k \}_{1 \leq k \leq \frac{q^n-1}{q-1}}.$
\item The only places of $F$ that are ramified in $E_{i,Kum}$ are $P_\infty$ and $\{P_\alpha\}_{s_{n,i}(\alpha)=0}$. Let $v_P = v_P(s_{n,i}(x))$.  Define $r_\alpha = \gcd([E_{i,Kum}:F],v_{P_\alpha}) > 0$ and $r_{\infty}= \gcd([E_{i,Kum}:F],v_{P_\infty}) > 0$. Then we have that
\begin{eqnarray*}
 e(P'_\alpha|P_\alpha) &=& \frac{[E_{i,Kum}:F]}{r_\alpha}, \\
 e(P'_\infty|P_\infty) &=& \frac{[E_{i,Kum}:F]}{r_{\infty}}.
\end{eqnarray*}
Since the extension is tame, the different exponents are given by
 $$d(P'|P) = e(P'|P) - 1,\; \forall P' \in \pe.$$
Also, $v_{P_0}(s_{n,i}(x))= \frac{q^i-1}{q-1}$ and $v_{P_\infty}(s_{n,i}(x))= q^{n-i}(\frac{q^i-1}{q-1})$.
\item All other rational places of ${\Bbb F}_{q^n}(x)$ split completely in $E_{i,Kum}.$
\end{enumerate}
\end{theorem}
{\bf Proof}. The proofs for (i) and (ii) will need standard results on Kummer extensions, cf. \cite{Sti1}, Ch. III.7.3. Also note that $r_\alpha, r_\infty < \frac{q^n-1}{q-1}$. For (iii) we note that from Lemma~\ref{lemma:fsphi-in-fq}, $s_{n,i}(\gamma) \in \fq, \forall \gamma \in \nq$, and therefore, it has $\frac{q^n-1}{q-1}$ pre-images under the norm map. \hfill $\Box$

\begin{lemma}
$$E_{i,Kum} \cong E_{n-i,Kum}.$$
\end{lemma}
{\bf Proof}. By making the substitution $y = xy'$, and then using Lemma~\ref{lemma:Sni-Sn(n-i)}. \hfill $\Box$

This immediately leads to the following corollary.
\begin{corollary} \label{corollary:isomorphic-to-trace-norm}
The extensions $E_{1,Kum}$ and $E_{n-1,Kum}$ of Theorem~$\ref{theorem:Kummer-symmetric}$ are both isomorphic to the trace-norm extension $E_n$ of Theorem~$\ref{theorem:Art-Sch-symmetric}$.
\end{corollary}

\begin{example} \rm The following are the symmetric Kummer extensions of ${\Bbb F}_{q^5}(x)$. $E_{1,Kum}$ is again the trace-norm extension.
$$
\begin{array}{|c|c|} \hline
E_{1,Kum} &  \mbox{(trace-norm)} y^{q^4 + q^3 + q^2 + q + 1}  = x^{q^4} +x^{q^3} + x^{q^2} + x^q + 1 \\ \hline
E_{2,Kum} &  y^{q^4 + q^3 + q^2 + q + 1} = x^{q^4+q^3} + x^{q^4 + q^2} + x^{q^4 + q} + x^{q^4 + 1} + \\
    & x^{q^3 + q^2} + x^{q^3 + q} + x^{q^3 +1} + x^{q^2+q} + x^{q^2 +1} + x^{q+1} \\ \hline
E_{3,Kum} &  y^{q^4 + q^3 + q^2 + q + 1} = x^{q^4+q^3+q^2} + x^{q^4+q^3 +q}+x^{q^4 + q^3 + 1} + x^{q^4 + q^2+q} \\
    & + x^{q^4 + q^2 + 1} + x^{q^4 + q + 1} + x^{q^3+q^2+q} + x^{q^3+q^2 +1} + x^{q^3+q+1} + x^{q^2+q+1} \\ \hline
E_{4,Kum} &  y^{q^4 + q^3 + q^2 + q + 1}  = x^{q^4+q^3+q^2+q} + x^{q^4+q^3+q^2+1} + x^{q^4+q^3+q+1} \\
          & + x^{q^4 + q^2 + q + 1} + x^{q^3+q^2+q+1} \\ \hline
\end{array}
$$
\end{example}

\begin{example} \rm (n=3, i=2)  Let $F=K(x)$, where $K = {\Bbb F}_{8}$. Let $E=F(y)$ where $y$ satisifes the equation
$$ y^7 = x^6 + x^5 + x^3 \; ( = s_{3,2}(x)). $$
$s_{3,2}(x)$ has three distinct zeros, other than $0$ itself, which has multiplicity three. $N(E) = 33$ and $g(E) = 9$. From Corollary~\ref{corollary:isomorphic-to-trace-norm}, we can see why these are the same values as the corresponding trace-norm extension.
\end{example}

There are other examples of Kummer extensions using $(n,q)$-elementary symmetric polynomials alone. An example from \cite{GarSti3} is given below.

\begin{example} \rm
Let $F = \fq(x)$, $q=p^e, e > 1, m = \frac{q-1}{p-1}$, and let $E=F(y)$ with
$$ y^m = (1+x)^m + 1.$$
Here $m$ rational places in \pf ramify, while all others split completely in $E$. Notice that $P_\infty$ splits completely in this extension, unlike the extensions of Theorem~\ref{theorem:Kummer-symmetric}. This construction, when iterated, gives an asymptotically good tower \cite{GarSti3}.
\end{example}

\section{Number of rational places versus genus} \label{section:rational-places-versus-genus}

In order to arrive at the Oesterle lower bound on genus $g$ of a function field for a specified field of constants \fq and specified number of rational places $N = L +1$, one must go through an algorithm of sorts, that is given below \cite{Sch2}.
\renewcommand{\labelenumi}{\rm (\arabic{enumi})}
\begin{enumerate}
\item Let $m$ be the unique integer for which
$$ {\sqrt{q}}^{\,m} < L \leq {\sqrt{q}}^{\,m+1}. $$
\item  Define
$$ u=\frac{{\sqrt{q}}^{\,m+1} - L}{L{\sqrt{q}}-{\sqrt{q}}^{\,m}} \in [0,1).$$
\item  Let $\theta_0$ be the unique solution of the trigonometric equation
$$ cos{\frac{m+1}{2}\theta} + u\, cos{\frac{m-1}{2}\theta} = 0 $$
in the interval $[\frac{\pi}{m+1}, \frac{\pi}{m}]$. \\
\item  Then we have that
$$g \geq \frac{(L-1)\sqrt{q}\, cos \theta_0 + q - L}{q+1-2\,\sqrt{q}\, cos \theta_0}.$$
\end{enumerate}
\renewcommand{\labelenumi}{\rm (\roman{enumi})}

\begin{example} \rm
If $q=8$ and $N=17$, then we have that $m=2$, $u=0.1779$ and $\theta_0 =1.1472$ yielding $g \geq 1.414$. Thus for a function field  over  ${\Bbb F}_{8}$ to have $17$ rational places, it must have genus at least 2.
\end{example}

\begin{example} \rm
If $q=8$, and $N=33$, then we have that $m=3$, $u=0.47$ and $\theta_0 = 0.88735$, yielding $g \geq 5.779$. Thus for a function field over  ${\Bbb F}_{8}$ to have $33$ rational places, it must have genus at least 6.
\end{example}

\begin{example} \rm
If $q=16$, and $N=129$, then we have that $m=3, u=0.2857$ and $
\theta_0=0.87752$, yielding $g \geq 17.88$. Thus for a function field over ${\Bbb F}_{16}$ to have 129 rational places, it must have genus at least 18.
\end{example}

\begin{example} \rm
If $q=27$, and $N=244$, then we have that $m=3, u=0.433$ and $
\theta_0=0.90754$, yielding $g \geq 25.16$. Thus for a function field over ${\Bbb F}_{27}$ to have 244 rational places, it must have genus at least 26.
\end{example}

Now we investigate the performance of symmetric function fields with respect to the known (Hasse-Weil, Oesterle) bounds on the number of rational places.

Theorem~\ref{theorem:Art-Sch-symmetric} describes, for a specific value of $n$, $n-1$ different symmetric extensions of the rational function field corresponding to  $i=2,3,\ldots,n$. The case of $i=n$ is the trace-norm extension. We now compare these $n-1$ different extensions for various values of $q$.
\begin{center}
\begin{tabular}{|c|c|c|c|c|c|c|c|} \hline
$n$ & $q$ & $N(E)$  & $g_n(E)$ & $g_{n-1}(E)$ & $g_{n-2}(E)$ & $g_{n-3}(E)$ & Oesterle  \\ \hline
3 & 2 & (17)  &     & $(2)^*$ & &               &$2^*$      \\
3 & 2 & 33    & 9   & $6^*$ &   &               &$6^*$      \\
3 & 4 & 1025  & 150     & 120   &   &                       & 74        \\
3 & 8 & 32769 & 2268    & 2016  &   &               & 903       \\
3 & 3 & 244   & 48      &  36   &   &                   & 26        \\
3 & 9 & 59050 & 3600    & 3240  &   &               & 1374      \\
3 & 5 & 3126  & 360 & 300   &   &               & 167       \\
3 & 7 & 16808 & 1344    & 1176  &   &               & 560       \\
3 & 11& 161052& 7920    & 7620  &   &               & 2808      \\
4 & 2 & 129   & 49  & 42    & 28    &               & 18        \\
4 & 4 & 16385 & 2646    & 2560  & 2016  &               & 667       \\
4 & 3 & 2188  & 507 & 468   & 351   &               & 152       \\
4 & 5 & 78126 & 9610    & 9300  & 7750  &               & 2071      \\
5 & 2 & 513   & 225 & 210   & 180   & 120               & 57        \\
\hline
\end{tabular}
\end{center}
{\bf Table 1}: {\footnotesize Number of rational points and genus for some members of the  families of symmetric function fields. The entry $g_j(E)$ denotes the genus of the extension described by Theorem~\ref{theorem:Art-Sch-symmetric} obtained for $i=j$. Note that $g_n(E)$ is the genus of the trace-norm extension. * indicates that the function field attains the Oesterle bound, and parentheses indicate a subfield of the symmetric function field.}

As  Theorem~\ref{theorem:Art-Sch-symmetric} tells us, the extension $E_2$ always has the lowest genus, while the extension $E_n$ (the trace-norm extension) has the highest genus. In other words, for purposes of the ratio of $N/g$, the trace-norm extension is actually the {\em worst} of the symmetric extensions, while the extension that uses the second $(n,q)$-symmetric polynomial is the {\em best}.

We now look more closely at the case of $n=3,\;i=2$, since this contains a maximal member in the the Oesterle sense. Specifically, we provide a uniformizing parameter at the place $P_\infty'$.
\begin{example} \rm \label{example:n=3,i=2,contd}
Let $F=K(x)$, where $K = {\Bbb F}_{q^3}$. Let $E=F(y)$ where $y$ satisfies the
equation
$$ y^{q^2} + y^q + y \;=\; x^{q^2+q} + x^{q^2+1} + x^{q+1}. $$
Except for $P_{\infty}$, which is totally ramified, all other rational places in $F$ split completely in $E$, giving a total of $q^5+1$ rational places. A uniformizing parameter for $P_\infty'$ can be obtained as follows. First we make the substitutions $x = X/Z \mbox{ and } y = Y/Z$ to get the following homogeneous equation:

$$Y^{q^2}Z^q + Y^qZ^{q^2} + YZ^{q^2+q-1} \;=\; X^{q^2+q}+X^{q^2+1}Z^{q-1}+X^{q+1}Z^{q^2-1}, $$
where, we look at $Y=1, Z = 0$ (\ie, $P'_\infty$), and we get
$$Z^q + Z^{q^2} + Z^{q^2+q-1} = X^{q^2 + q} + X^{q^2+1}Z^{q-1} + X^{q+1}Z^{q^2-1}. $$
Note that the valuations for $x,y,X,Y \mbox{ and } Z$ at $P'_\infty$ are $-q^2, -(q^2+q), q, 0 \mbox{ and } q^2+q$, respectively. Now, we expand $Z$ in a power series in $X$:
$$ Z = X^{q+1} + X^{2q} + \frac{X^{2q+1}}{\pi}, $$
where $\pi$ has valuation of the $q^{th}$ root of $X$, \ie, it is a uniformizing parameter for $P'_\infty$. Thus we get a uniformizing parameter $\pi$ for $P'_\infty$ given by:
$$ \pi = \frac{X^{2q+1}}{Z - X^{q+1} + X^{2q}} = \frac{x^{2q+1}}{y^{2q} - x^{q+1}y^q - x^{2q}y}. $$
Now, let $\sigma_\beta \in \gal(E/F)$ be such that $\sigma_\beta: y \rightarrow y + \beta$, where $\beta$ is a non-zero element of ${\Bbb F}_{q^3}$ whose trace in \fq is zero. Then,
$$ \sigma_\beta(\pi) - \pi = \frac{x^{2q+1}[\beta^{2q} + 2\beta^q y^q - \beta^q x^{q+1} - \beta x^{2q}]}{[(y+\beta)^{2q} - x^{q+1}(y+\beta)^q - x^{2q}(y+\beta)][y^{2q} - x^{q+1}y^q - x^{2q}y]}. $$
By noting that the valuation of both the terms in the denominator is the same, since they are conjugates, and then observing that the term in the numerator whose valuation will dominate is the last one, we get the valuation of $\sigma_\beta(\pi) - \pi$, which is independent of the choice of $\beta$:
$$ v_{P'_\infty}(\sigma_\beta(\pi) - \pi) = q^2 + 2. $$
This again gives us the sequence of ramification groups:
$$ G_0 = G_1 = \ldots = G_{q^2 +1}, $$
$$ G_{q^2 + 2} = \{0\}. $$
Thus the different exponent and the degree of the different for the extension $E/F$ are
$$ d(P_\infty'|P_\infty) =  \degdiff(E/F) = (q^2 + 2)(q^2 -1).  $$
This agrees with the results given in Theorem~\ref{theorem:Art-Sch-symmetric}.

\end{example}

We now wish to compare extensions of the rational function field having the typical (to our family of extensions) jumps in the sequence of ramification groups relative to the unique ramified place but having different degrees as extensions of the rational function field.
                                                                               \begin{theorem}
Let $F=K(x)$ and $K = \fq$. Let $E/F$ be an extension of function fields of degree $d$, where the sequence of ramification groups relative to the unique ramified place is of the form
$$ G_0=G_1=\ldots=G_k, $$
$$ G_{k+1} = \{0\}. $$
Further let all other places in $F$ split completely in $E$. Then the ratio of $N/g$ for $E$ decreases with increasing $d$.
\end{theorem}
{\bf Proof}. We have that
\begin{eqnarray*}
N &=& dq + 1, \\2g & = & 2 -2d + (d-1)(k+1) = (d-1)(k-1), \\
\frac{N}{2g} &=& \frac{dq + 1}{(d-1)(k-1)}.
\end{eqnarray*}
Then differentiating w.r.t $d$, we get that
$$ \left[\frac{N}{2g}\right]' = \frac{(q+1)(1-k)}{[(d-1)(k-1)]^2}. $$
Now notice that $k \geq 2$ for $g(E) > 0$. But in that case the numerator is always nonpositive. Hence the result. \hfill $\Box$

Thus, we get higher $N/g$ ratios by looking at subfields that have identical jumps in their sequence of ramification groups. The smaller the subfield, the higher is the $N/g$ ratio. In a sense, the rational function field is the limiting case of this behaviour, with $N/g = \infty$.

Thus ray class field extensions, which are maximal Abelian extensions with a certain property, have lower $N/g$ ratios than other Abelian extensions having identical jumps in their sequence of ramification groups.

\begin{example} \rm
 Consider the extension of Example~$\ref{example:n=3,i=2,contd}$. Let $q=2$ and  $G= \gal(E/F)$. Then the filtration of ramification groups at $P'_\infty$ is $G = G_0 = G_1 = \ldots = G_5$ and $G_6$ is trivial. In the ray class field constructions of function fields with many rational places \cite{Lau1}, there is a function field over ${\Bbb F}_8$ with an identical sequence of ramification groups. It has degree $8$, as an extension of the rational function field,  and genus $g=14$, for $N=65$ rational places. Thus,  the ratio  $N/g = 4.64$. The function field in our example has $N=33$ and $g=6$, giving a ratio of $N/g= 5.5$.
\end{example}

\section{Generalization of the Hermitian function field} \label{section:generalization-of-Hermitian}

 The simplest example of a symmetric function field is the Hermitian function field itself, described below.

Let $F={\Bbb F}_{q^2}(x)$ and $E=F(y)$ where $y$ satisfies the equation
$$ y^q + y \,=\, x^{q+1}. $$
Then $E$ is called the Hermitian function field. We have that $ g(E) = \frac{q(q-1)}{2}$, and $ N(E) = q^3 + 1$, which is  the maximum number allowed by the Hasse-Weil bound for this value of genus.  Furthermore, the Hermitian function field is the unique maximal function field over ${\Bbb F}_{q^2}$ of genus $g \geq \frac{q(q-1)}{2}$ \cite{FuhSti1}.  The Hermitian function field remains maximal for all values of $q$. Other symmetric function fields do not exhibit such a uniform performance. For instance, the function fields of  Example~\ref{example:n=3,i=2} and Example~\ref{example:n=3,i=2,subfield} attain the the Oesterle bound for $q=2$, while for higher values of $q$, they deviate considerably from the Oesterle bound. The other unique feature of the Hermitian function field is the extremely large size of its automorphism group \cite{Sti3}.

It has often been implicitly assumed that function field $E_n$ described by Theorem~\ref{theorem:Art-Sch-symmetric}, called the trace-norm function field, is the generalization of the Hermitian function field for $n \geq 3$. However, we argue below that it is the function field $E_2$ described by Theorem~\ref{theorem:Art-Sch-symmetric} that is more the generalization. It is clear that both coincide for $n=2$ and in this case, they are just the Hermitian function field. For $n \geq 3$, we must decide which is the more appropriate generalization of the Hermitian function field.  We now show the following similarities.

\vspace{0.2in}
\noindent{\bf A. Genus}\\
We have already observed that for  the $n-1$ families of symmetric extensions $E_i, 2  \leq i \leq n$ described by Theorem~\ref{theorem:Art-Sch-symmetric},
$$ g(E_2) < g(E_3) \ldots < g(E_n). $$
Thus, among these families $E_2$ has the lowest genus while the trace-norm function field has the highest.

\vspace{0.2in}
\noindent{\bf B. Automorphisms}\\
It is well known \cite{Sti2} that for every $\delta$ and $\tau$ in ${\Bbb F}_{q^2}$ that satisfy $\tau^q + \tau = \delta^{q+1}$, there is an automorphism $\sigma$ of the Hermitian function field  given by
\begin{eqnarray*}
\sigma(x)&=& x+ \tau, \\
\sigma(y)&=& y+ x \tau^q + \delta.
\end{eqnarray*}

We can state the general theorem about automorphisms of the symmetric function field $E=F(y)$ where $y$ satisfies (\ref{equation:Art-Sch-symmetric})  for $i=2$.
\begin{theorem} \label{theorem:generalized-Hermitian}
Let $F=K(x)$ where $K=\nq$. Let $E=F(y)$, where $y$ satisfies the equation
\begin{equation}
 y^{q^{n-1}} + y^{q^{n-2}} + \ldots+ y = s_{n,2}(x). \label{equation:generalized-Hermitian}
\end{equation}
Let  $(\delta, \tau)$ satisfy $\delta, \tau \in K$ and
$$  \delta^{q^{n-1}} + \delta^{q^{n-2}} + \ldots+ \delta = s_{n,2}(\tau). $$
Let $ m = \lfloor \frac{n-1}{2} \rfloor $. Then there exists an automorphism $\sigma$ of $E$ given by
$$ \begin{array}{lcl}
 \sigma(x) & = & x+\tau,  \\
 \sigma(y) & = & y + x\tau^q + x\tau^{q^2} + \ldots + x\tau^{q^{n-1}} + \delta.
\end{array} $$
The set of these automorphisms keeps $F$ fixed and forms a subgroup
$\Gamma$ of order $q^{2n-1}$ of the full automorphism group of the function field $E/K$. $\Gamma$ acts transitively on the set of finite rational places of $E$.
\end{theorem}
{\bf Proof}. On expanding $s_{n,2}(x + \tau)$, we observe that there are, apart from all the terms of $s_{n,2}(x)$ and $s_{n,2}(\tau)$, cross terms in $x$ and $\tau$ with all possible pairs of exponents $\{q^i,q^j\}$ for $x$ and $\tau$, where $ i \neq j,\;0 \leq i,j \leq n-1 $. Thus, there are $n(n-1)$ such cross terms, all distinct. On expanding
$[\sigma(y)]^{q^{n-1}} + [\sigma(y)]^{q^{n-2}} + \ldots + \sigma(y)$, we get
$n(n-1)$ cross terms in $x$ and $\tau$, all distinct, with the same exponent pairs as earlier. \hfill $\Box$

\vspace{0.2in}
\noindent{\bf C. Places of degree two}\\
Let us denote by $N_m(E)$ the number of rational places of degree $m$ in a function field $E$. Thus, $N_m$ is the number of Galois conjugacy classes of ``new" rational points of the curve over the extension of degree $m$ of the original field of definition (\ie, points that belong strictly to the extension, and not to the original field).  For the Hermitian function field over ${\Bbb F}_{q^2}$, $N_2 = 0$. This means that there are no ``new" rational places over  ${\Bbb F}_{q^4}$ --  the number of rational points over both the fields is $q^3 +1$. The same phenomenon occurs in each of the other ``Deligne-Lusztig'' curves, namely the Suzuki and the Ree curves.

\begin{lemma} \label{lemma:N_2=0}
Let the hypotheses be as in Theorem~$\ref{theorem:Art-Sch-symmetric}$.
Then  $N_2(E_i) \neq 0$ only if there exist ${\alpha, \beta}$, such that
$$ \beta^{q^{n-1}} + \beta^{q^{n-2}} + \ldots+ \beta = s_{n,i}(\alpha), $$
with  $\alpha \in {\Bbb F}_{q^{2n}} \setminus {\Bbb F}_{q^n}$ such that
$s_{n-1,i-1}(\alpha) \in {\Bbb F}_{q^n}$.
\end{lemma}
{\bf Proof}. Let $(x,y) = (\alpha,\beta)$ satisfy (\ref{equation:Art-Sch-symmetric}). Then, raising the equation to the $q^{th}$ power  and subtracting the original from it, we get
$$ \beta^{q^n} - \beta = [s_{n-1,i-1}(\alpha)]^q(\alpha^{q^n} - \alpha). $$
Now raising both sides to the power of $q^n$, we get that if $\alpha$ and $\beta$ are
elements of ${\Bbb F}_{q^{2n}}$, then
$$ \beta - \beta^{q^n} = [s_{n-1,i-1}(\alpha)]^{q^{n+1}}(\alpha - \alpha^{q^n}).$$
Then adding the two equations, we get that either $\alpha=\alpha^{q^n}$ or
$[s_{n-1,i-1}(\alpha)]^q=[s_{n-1,i-1}(\alpha)]^{q^{n+1}}]$. But since $\alpha \notin \nq$, and \nq is closed under taking $q^{th}$ powers, this gives us the result. \hfill $\Box$

\begin{example}
In the case of the Hermitian function field, since $s_{1,1}(x)=x$, Lemma~$\ref{lemma:N_2=0}$ precludes any possibility for $N_2 \neq 0$.
\end{example}

We now need another lemma.

\begin{lemma} \label{lemma:neq-2-3-6}
For $n \neq 3,4,6$, There exists no $\alpha \in {\Bbb F}_{q^{2n}} \setminus {\Bbb F}_{q^n}$ such that $s_{n-1,1}(\alpha) \in {\Bbb F}_{q^n}$.
\end{lemma}
{\bf Proof}.  Assume there is such an element $\alpha$. Then we have that
$$ \alpha^{q^{n-2}} + \alpha^{q^{n-3}} + \ldots + \alpha \in {\Bbb F}_{q^n},$$
which gives, on raising to the power of $q$
$$ \alpha^{q^{n-1}} + \alpha^{q^{n-2}} + \ldots + \alpha^q \in {\Bbb F}_{q^n}.$$
Subtracting the first equation from the second, we get
$$ \alpha^{q^{n-1}} - \alpha \in {\Bbb F}_{q^n}. $$
This implies that
$$ \alpha^{q^{2n-1}} - \alpha^{q^n} = \alpha^{q^{n-1}} - \alpha. $$
Now raising both sides to the power of $q^{n-1}$, we get
$$ \alpha^{q^{n-2}} - \alpha^{q^{2n-1}} = \alpha^{q^{2n-2}} - \alpha^{q^{n-1}}. $$
Adding this equation to the previous one, we get
$$ \alpha^{q^{n-2}} -  \alpha^{q^n} = \alpha^{q^{2n-2}} - \alpha. $$
This gives us
$$ \alpha^{q^n} - \alpha \in {\Bbb F}_{q^{n-2}}. $$

For this to hold, we would need ${\Bbb F}_{q^{n-2}}$ to be a subfield of ${\Bbb F}_{q^{2n}}$, or, $n-2\;|\;2n$. But for $n-2 > 0$,  this happens only for $n=3,4\mbox{ and } 6$. \hfill $\Box$

\begin{theorem} If $n \neq 3,4 \mbox{ or } 6$,
$$N_2(E_2) = 0.$$
\end{theorem}
{\bf Proof}. Follows from Lemma~\ref{lemma:N_2=0} and Lemma~\ref{lemma:neq-2-3-6}.  \hfill $\Box$

\vspace{0.2in}
We conclude this paper with an interesting insight into the curves corresponding to symmetric function fields. On any curve, for rational point $(\alpha, \beta)$ on the curve, the point's conjugate points  $({\alpha} ^{q^i}, {\beta}^{q^i}), 1 \leq i \leq n-1$, also lie on the curve. However, in general,  points
$({\alpha},{\beta}^{q^i}), 1 \leq i \leq n-1$, and $ ({\alpha}^{q^i}, {\beta}), 1 \leq i \leq n-1$, do not lie on the curve.  In the case of a curve
corresponding to a symmetric extension of the rational function field,  however, these  points too lie on the curve.  Thus we get many rational points ``for free". Similarly, we get many automorphisms also.

\section{Conclusion} \label{section:conclusion}
We have developed a theory that produces infinite families of function fields in which all finite rational places split completely. This theory works over all non-prime finite fields and thus explains, and does away with, the constraints of square cardinality on several existing examples and unifies them under one umbrella. The theory also provides a generalization of the Hermitian function field over non-square finite fields. Further research along the lines of generalizing the notions of symmetry developed in this paper have already been very fruitful leading to several new results (\cite{Deo2, Deo 3}) on the existence and explicit construction of function fields and towers with interesting properties.

\section*{Acknowledgements}
I would like to express my deep sense of gratitude to Prof. Dennis Estes, who supervised this work,  and tragically passed away just prior to its completion. Without his constant help, I could not have made any progress whatsoever. This work is dedicated to him.

I would also like to thank Joe Wetherell, who helped me immensely in completing this work following the demise of Prof. Dennis Estes.

\end{document}